\documentclass[Unicode]{cedram-alco}
\usepackage{amsmath,amsfonts,amssymb,amsthm,epsfig,enumitem,listings,mathrsfs,setspace}
\usepackage{color, bm}
\usepackage[vcentermath]{youngtab}
\usepackage{tikz}
\usetikzlibrary{matrix}
\tikzset{tab/.style={matrix of math nodes,column sep=-.4, row sep=-.4,text height=8pt,text width=8pt,align=center,inner sep=3}}
\tikzset{smalltab/.style={matrix of math nodes,column sep=-.4, row sep=-.4,text height=6pt,text width=6pt,align=center,inner sep=2.5}}

\usepackage{appendix}

\equalenv{Theorem}{thm}
\equalenv{Proposition}{prop} 
\equalenv{Lemma}{lemm}
\equalenv{Definition}{defi}
\equalenv{Example}{exam}
\equalenv{Remark}{rema}

\newcommand{\arxiv}[1]{\href{http://arxiv.org/abs/#1}{\tt arXiv:\nolinkurl{#1}}}

\newcommand{\wt}{{\rm wt}}

\newcommand{\Irr}{{\bm Irr}}
\newcommand{\bfv}{{\bm v}}
\newcommand{\Hom}{\text{Hom}}
\newcommand{\g}{\mathfrak{g}}
\newcommand{\im}{\operatorname{im}}

\newcommand{\vdim}{ \text{{\bf dim}}}

\begin{document}

\title[A quiver variety approach to root multiplicities]{A quiver variety approach to root multiplicities}

\author[\initial{P.} Tingley]{\firstname{Peter} \lastname{Tingley}}
\address{Loyola University Chicago \\ Dept. of Mathematics and Statistics \\
1032 W. Sheridan Rd. \\ Chicago, IL}
\email{ptingley@luc.edu}
\urladdr{http://webpages.math.luc.edu/~ptingley/}

\email{williams.colin.m@gmail.com}
 
\begin{abstract}
We present combinatorial upper bounds on dimensions of certain imaginary root spaces for symmetric Kac-Moody algebras. 
These come from the realization of the corresponding infinity-crystal using quiver varieties. The framework is general, but we only work out specifics in rank two. In that case we give explicit bounds. These turn out to be quite accurate, and in many cases exact,  even for some fairly large roots.
\end{abstract}


\thanks{This paper was inspired by a discussion with Alex Feingold at the satellite conference ``representation theory and related topics" to ICM 2014. I thank Alex for his insights, and the organizers for a  wonderful event. I also thank Alex for comments on a draft of this work. }
 
\maketitle

\section{Introduction}

Finite dimensional simple Lie algebras over ${\Bbb C}$ are often studied using the root space decomposition: the Lie algebra is the direct sum of the Cartan subalgebra and a number of 1-dimensional root spaces, which are the simultaneous eigenspaces of the Cartan subalgebra under the adjoint action. Kac-Moody algebras \cite{Kac:1990} are generally infinite dimensional but have similar behavior: a Cartan subalgebra is built into the definition, and the algebra is the direct sum of the Cartan and an infinite number of root spaces. The root spaces are no longer all 1-dimensional though. Their dimensions, the root multiplicities, usually grow quickly.

There has been considerable interest in these multiplicities. See \cite{CFL:2014} for a survey. Two exact methods are known to calculate them, both based on the Weyl-Kac denominator identity: a closed form formula due to Berman and Moody \cite{BM} and a recursive formula due to Peterson \cite{Pete:rec}. In special cases these have been further investigated and combinatorialized in \cite{FF, KLL, KM95}. None of this gives completely satisfactory information about asymptotics, and open questions remain. See \cite[Open Problems 2 and 3]{CFL:2014} and Frenkel's conjectural upper bound for hyperbolic cases \cite{Frenkel-conj} (although counter examples to this are known in $E_{10}$, see \cite{KMW} and also \cite[\S8]{CFL:2014}). 

Here we propose a new approach to root space multiplicities and their asymptotics. The method goes through the combinatorics of the crystal $B(-\infty)$ and its geometric realization using quiver varieties. 

The crystal $B(-\infty)$ is a set that parameterizes a basis for the upper triangular part of the universal enveloping algebra, along with some operators that approximate the Chevalley generators. It is usually defined algebraically, but it can be realized in a variety of ways. Here we use the realization from \cite{KS:1997} where the underlying set consists of irreducible components of the varieties of nilpotent representations of Lusztig's preprojective algebra from \cite[\S 12]{L91}. 

In \cite{BKT:2014} the category of representations of the preprojective algebra is studied using Harder-Narasimhan filtrations. Irreducible components where this filtration generically has only one step are called stable, and it is shown that the number of stable irreducible components of a given weight is a restricted Kostant partition function. If that weight is a root, and that root is not a multiple of a smaller root, it is exactly the root multiplicity. 
Our method is to calculate those root multiplicities by counting stable irreducible components. We translate this to a more combinatorial problem using Kashiwara's string data, which is a way of labeling each $b \in B(-\infty)$ with a word in the index set of the Dynkin diagram. Calculating the root multiplicities then amounts to counting words such that:
\begin{enumerate}

\item[(I1)] \label{eq:intropart1} the result is a valid string data, and

\item[(I2)] \label{eq:intropart2} the corresponding component is stable.

\end{enumerate}
At least in rank 2, these impose simple combinatorial restrictions, and counting words subject to those conditions gives an upper bound on the root multiplicity. 

For instance, consider the hyperbolic algebra with Cartan matrix
\begin{equation} \label{eq:ris3}
\left(
\begin{array}{rr}
 2 &-3\\
 -3&2
\end{array}
\right).
\end{equation}
This is sometimes called the Fibonacci algebra, see \cite{Feingold}. 
For an imaginary root of the form $m \alpha_1+n \alpha_0$ for $m$ and $n$ relatively prime, the root multiplicity is bounded by the number of rational Dyck paths from $(0,0)$ to $(n,m)$ such that, for any consecutive edges of length $a_k, a_{k+1}$, $\frac{a_{k+1} }{a_k} < \frac{3+\sqrt{5}}{2} \simeq 2.618$ (the square of the golden ratio). 

In Theorem \ref{th:good-cond} we give a generalization for this result to any symmetric, hyperbolic, rank two Kac-Moody algebra. In Theorem \ref{th:cond2} we refine this result by ruling out more paths that violate (I2). This tighter bound appears to be quite accurate. 

Root multiplicities are given by data satisfying (I1) and (I2) in any symmetric type, but translating to combinatorics is more difficult in higher rank, partly because, as in \cite{Littelmann:1998}, (I1) gets quite involved. It should also be possible to consider non-symmetric types, either by ``folding" as in e.g. \cite{Sav05}, by re-working things in terms of KLR algebras using \cite{BKM:2014,Kleshchev:2014,KMuth, McNamara:??,TW}, or by using variations of quiver varieties for non-symmetric types from \cite{GLS,NT}. 

Root multiplicities have previously been studied using crystals in \cite{KM95}, and using Dyck paths in \cite{KLL}. Nonetheless, we believe our results are quite different. In particular, our condition $\frac{a_{k+1} }{a_k} < \frac{3+\sqrt{5}}{2}$ has not previously appeared. 

This note  is organized as follows. In \S \ref{S:background} we review Kac-Moody algebras, crystals and quiver varieties. In \S\ref{S:RM-SCSD} we discuss stability conditions and string data, and state our main observation (see \S\ref{ss:key}). In \S\ref{sec:rank2} we work out details in rank 2, resulting in explicit upper bounds (Theorems \ref{th:good-cond} and \ref{th:cond2}). In \S\ref{sec:KLR} we discuss relationships with KLR algebras. In \S\ref{sec:heuristic} we discuss heuristics suggesting our bounds should be fairly tight. In an Appendix written with Colin Williams we present computational results.

\section{Background} \label{S:background}

\subsection{Kac-Moody algebras}
Fix a symmetric Cartan matrix $A$ with index set $I$. The Kac-Moody algebra $\mathfrak{g}$ is the Lie algebra generated by $\{ E_i, F_i, H_i : i \in I \}$ subject to the relations
\begin{itemize}

\item $[H_i, H_j]=0$,

\item $[H_i, E_j]= a_{ij} E_j$ and  $[H_i, F_j]= -a_{ij} F_j$,

\item $[E_i,F_j] = \begin{cases} H_i \quad i=j \\ 0 \quad \text{otherwise,} \end{cases}$

\item  For $i \neq j, \text{ad}_{E_i}^{-a_{ij}+1} E_j = 0 \quad \text{and} \quad  \text{ad}_{F_i}^{-a_{ij}+1} F_j = 0$. 
\end{itemize}
Here $ad$ is defined by $\text{ad}_X Y= [X,Y]$. As usual, let $\{ \alpha_i\}_{i \in I}$ be the simple roots, let $Q$ be their ${\Bbb Z}$-span,
and $Q_+$ their ${\Bbb Z}_{\geq 0}$ span. Then $\g$ is $Q$-graded, where, for each $i$, 
$$\deg E_i=-\deg F_i = \alpha_i, \ \ \deg H_i=0.$$

A non-zero $\beta \in Q$  is called a root if $\dim \g_\beta \neq 0$, in which case $m_\beta:= \dim \g_\beta $ is called the root multiplicity. All roots are either positive roots, meaning they are ${\Bbb Z}_{\geq 0}$ linear combinations of the simple roots $\alpha_i$, or negative roots, meaning the negatives of these. Let $\Delta$ denote the set of roots and $\Delta_+$ the positive roots.

There is an inner product on $Q$ defined by, for simple roots $\alpha_i,\alpha_j$, $\langle \alpha_i, \alpha_j \rangle = a_{ij}$. All roots $\beta$ have the property that either $\langle \beta, \beta \rangle=2$, in which case $\beta$ is called a real root, or $\langle \beta, \beta \rangle \leq 0$, in which case $\beta$ is called an imaginary root.

Let $U(\mathfrak{g})$ be the universal enveloping algebra of $\mathfrak{g}$. As a vector space, 
$$U(\g)= U^-(\g) \otimes U^0(\g) \otimes U^+(\g),$$
where $U^-, U^0, U^+$ are the subalgebras generated by the $F_i$, the $H_i,$ and the $ E_i$ respectively. The graded dimension of $U^+$ is 
$$\dim U^+ = \prod_{\beta \in \Delta_+} \left( \frac{1}{1-e^\beta} \right)^{m_\beta}.$$
That is, the dimension of the $\gamma$-weight space of $U^+(\g)$ is the number of Kostant partitions of $\gamma$, meaning the number of ways to write $\gamma$ as a sum of positive roots, taking into account multiplicities.

\subsection{The crystal $B(-\infty)$}
For any symmetrizable Kac-Moody algebra, the crystal $B(-\infty)$ is a set along with operators $e_{i}, f_{i} \colon B(-\infty) \rightarrow B(-\infty) \cup \{ 0 \}$ for each $i \in I$, which satisfy various axioms. Roughly, $B(-\infty)$ parameterizes a basis for $U^+(\mathfrak{g})$, and the $e_i, f_i$ are related to the Chevalley generators $E_i, F_i$. There is a weight function $\wt: B(-\infty) \rightarrow Q$, and the number of elements of a given weight $\gamma$ is the dimension of the $\gamma$ weight space in $U^+(\mathfrak{g})$.
See \cite{Kashiwara:1995} or \cite{HK:2002} (these sources consider $B(\infty)$, which is related to $B(-\infty)$ by Cartan involution). 
Here we only need the realization of $B(-\infty)$ from \cite{KS:1997}, which is explained below.

\subsection{Quiver varieties } \label{ss:QV}

Fix a graph $G$ with vertex set $I$ and edge set $E$. Let $\overline Q$ be the corresponding double quiver, which is the directed graph with vertex set $I$ and arrow set $A$ consisting of two arrows for each edge $e \in E$, one in each direction. For each arrow $a \in A$, let $s(a)$ be the source and $t(a)$ be the target, meaning $a$ points from $s(a)$ to $t(a)$. 
\begin{Definition}
The path algebra ${\Bbb C} \overline Q$ is the algebra over ${\Bbb C}$ with basis consisting of all paths in $G$ (sequences of arrows $a_k \cdots a_1$ with $t(a_i)=s(a_{i+1})$), plus the lazy paths $e_i$ at each vertex) and with multiplication given by 
$$(b_{k} \cdots b_1) (a_j \cdots a_1) = 
\begin{cases}
b_k \cdots b_1 a_j \cdots a_1 \quad \quad t(a_j)=s(b_1) \\
0 \quad \quad \quad \text{otherwise}.
\end{cases}
$$
\end{Definition}
Choose a subset $\Omega$ of $A$ where each edge appears in exactly one direction (this is sometimes called an orientation of $G$). Define $\varepsilon(a)= 1$ if $a \in \Omega$ and $-1$ otherwise. For any arrow $a$, let $\bar a$ denote the reverse of $a$. 

\begin{Definition}
The preprojective algebra $\Lambda$ is the quotient of ${\Bbb C}\overline Q$ by the ideal generated by
$$\rho =   \sum_{a \in A} \varepsilon(a) \bar a a.$$
\end{Definition}

\begin{Definition}
For any $I$-graded vector space $V = \oplus_{i \in I} V_i$, let $\Lambda(V)$ be the variety of actions of $\Lambda$ on $V$ where the lazy path $e_i$ at $i$ acts as projection onto $V_i$, and which are nilpotent in the sense that all paths of length at least $\dim V$ act as $0$.  
\end{Definition}
We call $\Lambda(V)$ the quiver variety. It is also sometimes called Lusztig's nilpotent variety, to distinguish it from other quiver varieties in the literature. 
A point $x \in \Lambda(V)$ is an algebra homomorphism $p \rightarrow x_p$ from $\Lambda$ to $\text{End}(V)$ such that $x_{e_i}$ is the projection $\pi_i$ onto $V_i$, and such that $x_p=0$ for all sufficiently long paths $p$.
Notice that $x$ is determined by $\{x_a\}_{a \in A}$, and each $x_a$ is in $ \Hom (V_{s(a)}, V_{t(a)})$. In this way $\Lambda(V)$ is a subvariety of $\oplus_{a\in A} \Hom (V_{s(a)}, V_{t(a)})$.

Up to isomorphism, $\Lambda(V)$ depends only on ${\bfv} = \vdim V = (\dim V_i)_{i \in U}$. 
Identify $\bfv$ with the point $\gamma= \sum_{i\in I} v_i \alpha_i$ in the root lattice, and denote $\Lambda(V)$ by $\Lambda(\gamma)$. Let $\Irr \Lambda (\gamma)$ denote the set of irreducible components of $\Lambda(\gamma)$.

Associate to each graph a symmetric Cartan matrix whose index set is the set of vertices, and where, for $i \neq j$, $-a_{ij}$ is the number of edges connecting $i$ and $j$. 
The following is due to Kashiwara and Saito \cite{KS:1997}, and can be found in the current form in \cite{NT}.

\begin{Theorem} \label{th:QVrealization}
The crystal $B(-\infty)$ is naturally indexed by 
$\coprod_{\gamma \in Q_+} \Irr \Lambda (\gamma)$. The operation $f_i^{\text{max}}$ which applies the crystal operator $f_i$ as many times as possible acts on $X \in \Irr \Lambda (\gamma)$ as follows:
Fix $T \in X$. Let $\text{Soc}_i(T)$ be the intersection of the socle of $T$  with $T_i$ and set $\gamma'= \gamma-\dim \text{Soc}_i(T) \alpha_i$. Generically $T/{\text{Soc}_i}(T)$ is isomorphic to a point in a unique $Y \in \Irr \Lambda(\gamma')$, and $f_i^\text{max} X =Y$.   \qed
\end{Theorem}

\begin{Example} \label{ex:beginning}
Here the most important example is the graph
\begin{center}
\begin{tikzpicture}
\node [draw, circle, fill=black] at (0,0) {.};
\node [draw, circle, fill=black] at (3,0) {.};

\draw[line width = 0.02cm] (0,-0.14)--(3,-0.14);
\draw[line width = 0.02cm] (0,0)--(3,0);
\draw[line width = 0.02cm] (0,0.14)--(3,0.14);
\end{tikzpicture}
\end{center}
corresponding to the ``Fibonacci" algebra with Cartan matrix
$$
\left(
\begin{array}{rr}
2 & -3 \\
-3 & 2
\end{array}
\right).
$$
Oriente left to right, which is to say choose $\Omega$ to consist of all arrows pointing left to right. Consider $V$ with dimension vector $\bfv=(3,4)$. An action of ${\Bbb C}\overline Q$ on $V$ is defined by three maps $x^1,x^2,x^3 : {\Bbb C}^3 \rightarrow {\Bbb C}^4$, one for each arrow pointing left to right, and three maps  $y^1,y^2,y^3 : {\Bbb C}^4 \rightarrow {\Bbb C}^3$, corresponding to the reverse arrows. So the representation variety of ${\Bbb C}\overline Q$ in this dimension is isomorphic to ${\Bbb C}^{6 \times 12}$. $\Lambda(v)$ is the sub-variety cut out by the condition that all paths of length 7 act as 0 and the equations
$$x^1 y^1+x^2 y^2+x^3 y^3=0, \quad y^1 x^1+y^2 x^2+y^3 x^3=0,$$
where the left equation is in $\text{End}{\Bbb C}^4$ and the right in  $\text{End} {\Bbb C}^3$.   
\end{Example}

\section{Root multiplicities from stability conditions and string data} \label{S:RM-SCSD}

\subsection{Stability conditions}

The following loosely follows \cite{BKT:2014}, and also draws on notation from \cite{TW}. 

\begin{Definition}
A {\bf charge} $c$ is a linear function $c: \mathfrak{h}^* \rightarrow {\Bbb C}$ such that the images of all simple roots (and hence all positive roots) are in the upper half plane. 
\end{Definition}

For a fixed charge $c$, any representation $T$ of $\Lambda$ has a unique filtration 
$$\emptyset = T_0 \subset T_1 \subset \cdots \subset T_k= T$$
where the sub-quotients $Q_i=T_i/T_{i-1}$ satisfy
\begin{enumerate}
\item[(HN1)] $Q_i$ has no submodule $S$ with $\arg(c(\vdim S))< \arg(c( \vdim Q_i))$,

\item [(HN2)] $\arg(c (\vdim Q_1)) < \arg(c (\vdim Q_2)) < \cdots < \arg(c (\vdim Q_k))$.

\end{enumerate}
Here $\arg$ is the angle in the complex plane. This is a special case of a Harder-Narasimhan filtration as in e.g. \cite{Rudakov97}, so we call it the HN filtration. The following follows by applying \cite[Theorem 4.4]{BKT:2014} repeatedly.

\begin{Theorem} \label{th:ss}
Fix an irreducible component $X$ of $\Lambda(V)$. For generic $T \in X$ each sub-quotient $T_j/T_{j-1}$ lies in a unique irreducible component $X_j$ of $\Lambda(\mathring{Q}_j)$. Furthermore, $X_j$ is generically constant, meaning it is constant on an open dense subset of $X$. Here $\mathring{Q}_j$ means the vector space ${Q}_j$ with the action of $\Lambda$ forgotten.  
\qed
\end{Theorem}
\noindent We call $X$ {\it stable} if the HN filtration for generic $T \in X$ has one step. This implies
\begin{itemize}
\item[(S)] For any submodule $S \subset T$, $\arg(c(\vdim S)) \geq \arg(c (\vdim T))$.
\end{itemize}
By Theorem \ref{th:ss} each $X \in \Irr \Lambda(\bfv)$ has a unique {\it stable decomposition } $(X_1, \ldots, X_k)$, where $X_i$ is stable and for generic $T \in X$, $Q_i$ is in $X_i$. 

Fix a stability condition $c$ so that, for any root $\beta$, if $\arg c(\alpha)=\arg c(\beta)$ then $\beta$ and $\alpha$ are parallel. The following can be extracted from \cite{BKT:2014}, and the proof below can be found in \cite[Corollary 2.12]{TW} in a somewhat different context.

\begin{Theorem} \label{th:BigK}
For any $\gamma \in Q_+$, the number of stable irreducible components of $\Lambda(\gamma)$ is the sum over all ways of writing $\gamma=v_1\beta_1+\cdots+v_n \beta_n$ as a sum of parallel roots $\beta_k$
of the product $m_{\beta_1} \cdots m_{\beta_n}$ of the corresponding root multiplicities. In particular, if $\gamma$ is not parallel to any smaller weight, the number of stable irreducible components is exactly $m_\gamma$. 
\end{Theorem}

\begin{proof}
If $\nu$ is a simple root the result is trivial.
  Proceed by induction on
  $\rho^\vee(\nu)$. 
 \begin{equation*} \label{eq:cp}
|\Irr \Lambda(\nu)| = \dim U^+_q(\g)_\nu =\sum_{\substack{\nu=\beta_1+\cdots+\beta_n}} \prod_{i=1}^n m_{\beta_i},
  \end{equation*}
 where the sum if over all tuples of positive roots whose sum if $\nu$. Inductively, components that have a semi-stable decomposition with at least two parts account for all the terms where the $c(\beta_j)$ do not all have the same argument. Thus the remaining terms where all the $\arg c(\beta_j)$ are equal, and hence all the $\beta_j$ are parallel, give the number of stable components.
\end{proof}

\subsection{String data} \label{ss:stringdata}
The following parameterization of $B(-\infty)$ is based on work of Berenstein and Zelevinsky \cite{BZ93} and Kashiwara 
\cite[\S 8.2]{Kashiwara:1995}, although conventions differ between those papers. Similar constructions were also studied by Nakashima and Zelevinsky \cite{NZ97} and Littelmann \cite{Littelmann:1998}. See e.g. \cite{FO17} for a discussion of how the different conventions relate. Here we follow \cite{BZ93, Littelmann:1998}.

Choose a sequence 
$i_1, i_2, i_3 \ldots$
of nodes in the Dynkin diagram with each appearing infinitely many times. 
The {\bf string data} $(a_1,a_2, \ldots)$ of $b \in B(-\infty)$ is defined by
\[
\begin{aligned}
& a_1= \max \{ n : f_{i_1}^n b \neq 0 \}, \\
&  a_2= \max \{ n : f_{i_2}^n f_{i_1}^{a_1} b \neq 0\},
\end{aligned}
\]
and so on. 
We record the string data as a word in the letters $I$ consisting of $a_1$ $i_1$'s, followed by $a_2$ $i_2$'s, and so on.  Sometimes we write this as
$$i_1^{a_1} i_2 ^{a_2} \cdots i_k ^{a_k}.$$

Indexing $B(-\infty)$ by $\sqcup \Irr \Lambda(\bfv)$, Theorem \ref{th:QVrealization} shows that the string data of $X \in\Irr \Lambda(\bfv)$ records the dimensions of the (graded) socle filtration of a generic $T \in X$:
$$
\begin{aligned}
& a_1 = \dim \Hom({\Bbb C}_{i_1}, T), \\
& a_2= \dim \Hom({\Bbb C}_{i_2}, T/\text{Soc}_{i_1}(T)) ,\\
\end{aligned}
$$
and so on, where ${\Bbb C}_i$ is the one dimensional simple module in degree $i$. 
In this way the notion of string data extends to all nilpotent $\Lambda$-modules. 
The string data in the crystal sense is the generic value of the string data on an irreducible component.

\subsection{Key Observation} \label{ss:key}
For a root $\beta$ which is not a multiple of a smaller root, Theorem \ref{th:BigK} shows that the root multiplicity $m_\beta$ is the number of string data which correspond to stable components, or equivalently the number of words satisfying (I1) and (I2). Describing these words combinatorially seems hard but, at least in rank two, we find a somewhat bigger set of words which can be understood combinatorially. The size of that set gives an upper bound on the root multiplicities. 

\section{Rank 2} \label{sec:rank2}
We now restrict to considering a Kac-Moody algebra with Cartan matrix
$$
\left(
\begin{array}{rr}
 2 &-r\\
 -r&2
\end{array}
\right)
$$
for $r \geq 3$, with $I =\{ 0, 1 \}$. Fix a charge $c$ with $\arg c(\alpha_0) < \arg c(\alpha_1)$, and take string data using the sequence $1,0,1,0,\ldots$.
Fix an imaginary root $\beta =m \alpha_1+n \alpha_0$ with $\text{gcd}(m,n)=1$, so $m_\beta$ is  the number of words in $\{0,1\}$ satisfying (I1) and (I2). 

Fix a $\Lambda$ module $T = T_0 \oplus T_1$ of dimension $\gamma$, and let
$$
\begin{aligned}
& U_1 = \text{Soc}_1(T) = \text{Soc}(V) \cap T_1 \\
& U_2 = U_1 \cup  \{ v \in T_0 : x_a(v) \subset U_1 \text{ for all } a \in A\}  
\end{aligned}$$
and so on be the graded socle filtration of $T$. As in \S\ref{ss:stringdata} the string data is $a_k = \dim U_k/U_{k-1}$. 
We want to count stable irreducible components. Since stability is an open condition we should count string data of weight $\beta$ such that there exists a stable module $T$ with that data. 

\subsection{Dyck path condition}
Fix $(a_1,a_2, \ldots)$, and assume there exists a stable $T$ with this string data. 
Each $U_{2k}$ is a submodule, so the stability condition (S) implies that, for all $k$, $$\frac{a_2 +a_4 \cdots+a_{2k}}{a_1 + a_3 \cdots+a_{2k-1}} \leq  \frac{n}{m}.$$ Draw the data as a path in the plane with $a_1$ steps up, then $a_2$ to the right, etc. This says that $(a_1,a_2, \ldots)$ is a rational Dyck path. That is, it does not go below the diagonal as shown:

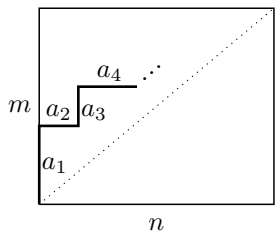
\begin{center}
\begin{tikzpicture}[scale=0.52]

\draw[line width= 0.02cm] (0,0)--(6,0)--(6,5)--(0,5)--(0,0);
\draw[line width= 0.035cm] (0,0)--(0,2)--(1,2)--(1,3)--(2.5,3);
\draw[line width= 0.015cm, dotted] (0,0)--(6,5);

\node at (2.7,3.2) {.};
\node at (2.88,3.35) {.};
\node at (3.06,3.5) {.};

\node at (0.4,1) {$a_1$};
\node at (0.5,2.35) {$a_2$};
\node at (1.4,2.35) {$a_3$};
\node at (1.8,3.35) {$a_4$};

\node at (3,-0.5) {$n$};
\node at (-0.5,2.5) {$m$};
\end{tikzpicture}
\end{center}

\noindent A simple argument dating to at least Grossman \cite{Grossman} shows that, since $m$ and $n$ are relatively prime, the number of such paths is  
\begin{equation} \label{eq:numberDyck}
\frac{1}{m+n} \left( \begin{array}{c} m+n \\ n \end{array} \right).
\end{equation}

\subsection{Condition on consecutive edge lengths} \label{condse}
We must also restrict to valid string data. String data have been characterized by Littelmann in rank 2:

\begin{Theorem} \cite[Proposition 2.1]{Littelmann:1998} \label{th:litth}
 $a_1,a_2,a_3, \ldots$ is a string data of some $b \in B(-\infty)$ if and only if
\begin{equation}
a_3 \alpha_0 \leq a_2 s_0 \alpha_1, \;\; a_4 s_0 \alpha_1 \leq a_3 s_0 s_1 \alpha_0, \;\; a_5 s_0 s_1 \alpha_0 \leq a_4s_0 s_1 s_0 \alpha_1, \ldots
\end{equation}
\end{Theorem}

As we will see, this condition becomes simpler for stable irreducible components.

\begin{Lemma} \label{lemma-1down}
Assume that for some $1\leq j<k$ with $k-j$ even, the sub-quotient $Q= U_k/U_{j}$ of $T$ has a submodule of dimension $\nu= a \alpha_{j}+b \alpha_{j+1}$. Then $Q'= U_{k-1}/U_{j-1}$ has a submodule of dimension $b \alpha_{j+1}+c \alpha_j$ for some $c \leq rb-a$. If $j\geq 2$ then $c \geq \frac{b}{r}$.
\end{Lemma}

\begin{proof}
For simplicity assume $j$ is odd, so $U_j/U_{j-1}$ and $U_k/U_{k-1}$ are in degree $1$. The case when $j$ is even follows by reversing the roles of $T_0$ and $T_1$. 

Clearly it suffices to assume the submodule is all of $Q$.  
Consider the map 
$$M=\bigoplus_{a \in A: s(a)=1} x_a: U_k \cap T_1 /U_{j} \cap T_1 \rightarrow (U_{k-1} \cap T_0 / U_{j-1} \cap T_0)^{\oplus r}.$$ 
This must be injective, as otherwise the kernel would be further down in the socle filtration, contradicting the definition of the submodules $U_i$ (in particular $b \neq 0$). Now consider the map
$$\bar M = \sum_{a \in A: s(a)=0} x_{\bar a} : (U_{k-1} \cap T_0/ U_{j-1} \cap T_0)^{\oplus r} \rightarrow  U_{k-2} \cap T_1 / U_{j-2} \cap T_1. $$
The preprojective relation implies this descends to a map 
$$ (U_{k-1} \cap T_0 / U_{j-1} \cap T_0)^{\oplus r}/ \im M \rightarrow U_{k-2} \cap T_1 / U_{j-2} \cap T_1,$$
and its image can have dimension at most the dimension of the domain, $ r b - a.$
Thus 
$$T= (U_{k-1} \cap T_0/U_{j-1} \cap T_0) \oplus \im \bar M$$
is a submodule of $Q'$ of dimension $b\alpha_{j+1}+c\alpha_j$ for $c\leq rb-a$ (since $j$ is odd, $U_{j-2} \cap T_1 = U_{j-1} \cap T_1$). 

If $j \geq 2$, by the definition of the socle filtration, the map
$$\bigoplus_{a: s(a)=0} x_a: T_0 \rightarrow T_1^{\oplus r}$$
is injective. This implies $c\geq \frac{b}{r}$.  
\end{proof}

\begin{Theorem} \label{th:good-cond}
Assume $a_1, a_2, \ldots a_N$ is the string data of a stable irreducible component. Let $\sum_{k=1}^N a_k \alpha_k = n\alpha_0+m\alpha_1$, and assume $\langle n\alpha_0+m\alpha_1, n\alpha_0+m\alpha_1  \rangle<0$. Here $\alpha_k$ means $\alpha_0$ if $k$ is even and $\alpha_1$ if $k$ is odd. Then, for all $k \geq 1$,
\begin{equation}  \label{eq:hyp-cond1}
\frac{a_{k+1}}{a_k} \leq \frac{\sqrt{r^2-4}+r}{2}.
\end{equation} 
In particular, if $m$ and $n$ are relatively prime and $\beta= n \alpha_0+ m \alpha_1$ is an imaginary root, then 
the number of rational Dyck paths from $(0,0)$ to $(n,m)$ satisfying \eqref{eq:hyp-cond1} for all $k$ is an upper bound for $m_\beta$.
\end{Theorem}

\begin{proof}
By evaluating the inner product, \eqref{eq:hyp-cond1} is equivalent to
\begin{equation} \label{eq:hrk}
a_{k+1} \leq a_k \quad \text{ or }\quad  \langle a_k \alpha_k + a_{k+1} \alpha_{k+1}, a_k \alpha_k + a_{k+1} \alpha_{k+1}\rangle \leq0.
\end{equation}
Fix data violating \eqref{eq:hrk} for some $k$, and assume a module $T$ has that string data. It suffices to show that $T$ has a submodule violating stability. 

Proceed by induction on $k$. The case $k=1$ is clear since if $a_2, a_1$ fail to satisfy \eqref{eq:hyp-cond1}, then
$a_2>a_1$ and $\langle a_1 \alpha_1+a_2 \alpha_0, a_1 \alpha_1+a_2 \alpha_0 \rangle >0$, which implies 
 $\displaystyle \frac{a_2}{a_1} > \frac{n}{m}$. 

Assume $a_{k+1}, a_k$  violates the condition for some $k >1$. By Lemma \ref{lemma-1down}, $T$ has a submodule 
whose string data violates the conditions for $k-1$, since replacing $a_k \alpha_k+ a_{k+1} \alpha_{k+1}$ with  $a_k \alpha_k+ (ra_k-a_{k+1}) \alpha_{k+1}$ preserves the condition $\langle \nu, \nu \rangle >0$ (it is reflection), and lowering the smaller coefficient also preserves this condition. By induction there is a submodule violating stability.
\end{proof}

\begin{Remark}
Any Dyck path satisfying \eqref{eq:hyp-cond1} also satisfies the conditions of Theorem \ref{th:litth}: The roots involved are 
$$\alpha_0, 3\alpha_0+\alpha_1, 8\alpha_0+3 \alpha_1, 21 \alpha_0+8 \alpha_1, 55 \alpha_0+21\alpha_1, \cdots$$
The coefficients are always Fibonacci numbers $F_{2n}$ and $F_{2n-2}$, and the condition follows from the fact that
$\frac{F_{2n}}{F_{2n-2}}$ is bounded below by the square of the golden ratio. 
\end{Remark}

\begin{Example} 
Consider $r=3$ and $\beta=4 \alpha_0 + 3 \alpha_1$. There are 
$\displaystyle \left( \begin{array}{c}  7 \\ 3 \end{array} \right)=35$ possible words. By Theorem \ref{th:litth}, all except
$$
1000011, 1010001, 1101000
$$
are string data for $1,0,1,0 \ldots$. For instance, $1010001$ violates the conditions  
because
$$3 (3 \alpha_0+\alpha_1) \not \leq 1 (8 \alpha_0+3 \alpha_1).$$
This correctly predicts $\dim U^+(\mathfrak g)_{3 \alpha_0 + 4 \alpha_1}=32$. There are only five rational Dyck paths, so five candidates for stable components:
$$
1110000, \quad
1101000, \quad 
1100100, \quad
1011000, \quad
1010100.$$
The path $1101000$ violates condition \eqref{eq:hyp-cond1} so does not correspond to a stable component. In fact, as above, this path does not even correspond to valid string data. 
\end{Example}

\begin{Example} \label{ex:part2}
Continue with $r=3$, but now consider $\beta=3 \alpha_0 + 4 \alpha_1$. The root multiplicity is still 4, and there are still 5 Dyck paths:
$$1111000, \quad 1110100, \quad 1110010, \quad 1101100, \quad
1101010.$$
This time they all satisfy \eqref{eq:hyp-cond1}, so our upper bound is off by one. 
The unstable component corresponds to $1101100$ since, by Lemma \ref{lemma-1down}, the sub-quotient corresponding to the subword $011$ implies the existence a submodule of the form $10$ or $0$, and that violates stability.
\end{Example}

In general our estimates are better in cases $m \alpha_0+n \alpha_1$ with $m$ slightly greater than n. The difference is largely explained by looking at the ends of the paths. If $m>n$, then Dyck paths must end $*00$, whereas if $m<n$ then they must only end $*0$. The last $0$s in a Dyck path can only cause a violation of \eqref{eq:hyp-cond1} if there are at least 3 of them, and this is more likely if $m>n$. 

\subsection{More restrictions} \label{ss:higher}
We now discuss a refinement to Theorem \ref{th:good-cond} giving a tighter upper bound. We start with some examples.

\begin{Example} For $r=3$ and the root  $8 \alpha_0+7\alpha_1$ the following satisfies \eqref{eq:hyp-cond1} 
but corresponds to a non-stable component: 
\begin{equation}
 1^2 0^2 1^5  0^6.
 \end{equation}
To see why this component is not stable, notice that, by 
Lemma \ref{lemma-1down}, the
sub-quotient $Q$ corresponding to the middle $0^21^5$ implies the existence of a submodule with string data $1^a0^2$ for $a=0$ or $1$. This 
violates stability.  
This path violates the conditions in Theorem \ref{th:cond2} below for $x=y=1$.
\end{Example}

\begin{Example} \label{ex:type2}
A similar problem can occur using two consecutive steps, and can occur further into the word. For instance, for the root
$16 \alpha_0+15\alpha_1$,
\begin{equation}
 {\color{red} 1^30^2}  1^2 {\color{red} 0^2} {\color{blue} 1^5}  {\color{red} 0^2} {\color{blue} 1^5}  0^{10}
\end{equation}
fails to be stable by looking at the submodule generated by the red numbers. 
This path violates the conditions in Theorem \ref{th:cond2} below for $x=2, y=3$.
\end{Example}

\begin{Theorem} \label{th:cond2}  
Fix a module $T$ with string data $a_1, a_2, \ldots, a_{2k}$. Let $m= a_1+a_3+\ldots+a_{2k-1}, n= a_2+a_4+\cdots+a_{2k}$. 
If $T$ is stable then, for all $1 \leq x\leq y< k$, 
\begin{equation} \label{type-2-all}
\frac{a_2+\cdots+ a_{2y}}{a_1+ \cdots+ a_{2x-3}+r(a_{2x} +\cdots+ a_{2y})-a_{2x+1} - \cdots - a_{2y+1}}
\leq
\frac{n}{m}.
\end{equation}
Here in each $\cdots$ the indices increase by 2 at a time. In particular, the number of Dyck paths satisfying \eqref{type-2-all} along with \eqref{eq:hyp-cond1} is a tighter upper bound on $m_\beta$. 
\end{Theorem}

\begin{proof}
Assume a module $T$ has string data where \eqref{type-2-all} is violated for some $x,y$. It suffices to show that $T$ is not stable.  
Applying Lemma \ref{lemma-1down} to the sub-quotient $U_{2y+1}/U_{2x-1}$ implies that 
$U_{2y}/U_{2x-2}$ has a submodule of dimension
$$
(a_{2x} + \cdots + a_{2y}) \alpha_0+ (r(a_{2x} + \cdots a_{2y}) - a_{2x+1}- \cdots - a_{2y+1}-k) \alpha_1$$
for some $k\geq0$. 
Taking this along with all of $U_{2x-2}$ gives a submodule of dimension
$$(a_2+ \cdots + a_{2y}) \alpha_0+ (a_1+a_3+\cdots + a_{2x-3}+ r(a_{2x} + \cdots a_{2y}) - a_{2x+1}- \cdots - a_{2y+1}-k) \alpha_1. $$
Since \eqref{type-2-all} is violated this submodule violates stability. 
\end{proof}

\begin{Example} \label{ex:type3}
For $r=3$ and  $16 \alpha_0+15 \alpha_1$, the calculations in the appendix show that there is exactly one word that satisfies both Theorem \ref{th:good-cond} and Theorem \ref{th:cond2} but does not correspond to a stable component. It is
$$1^{10}0^3 1^50^{13}. $$
To see that this does not correspond to a stable component, apply Lemma \ref{lemma-1down} to $1^50^{13}$ to obtain a sub-quotient of $U_3/U_1$ of the form $0^21^5$ ($2$ is the only integer between $\frac{5}{3}$ and $3*5-13$), and hence a submodule of the form $1^{10}0^21^5$. Applying Lemma \ref{lemma-1down} again gives a submodule $1 0^2$ or just $0^2$, violating stability.
\end{Example}

There are even stranger examples, but they seem to be exceedingly rare. 

\section{Relation with KLR algebras} \label{sec:KLR}

There is a version of this story using KLR algebras. By \cite[Corollary 2.12]{TW}, the root multiplicity of $m \alpha_0+n \alpha_1$ for $m,n$ relatively prime is the number of cuspidal representations of weight $m \alpha_0+n \alpha_1$ for a KLR algebra. These can be indexed by good Lyndon words (see \cite{HMM,Kleshchev-Ram}), so the root multiplicity is the number of such words. The problem is there is no nice combinatorial description of good Lyndon words. 

What we do here corresponds to instead labeling the cuspidal modules by their string data. This has a big advantage: at least in rank two, describing which words are string data is relatively easy. The string data itself is a word in the character of the module, so if the string data is not a Dyck path then the corresponding module is not cuspidal. But a string data which is a Dyck path can nonetheless correspond to a module which is not cuspidal, so we get an overestimate. We can make progress to correct the over-counting, and end up with estimates that seem quite good. But finding tractable conditions that exactly characterize string data of cuspidal modules seems tricky.

\section{Heuristics}  \label{sec:heuristic}

\subsection{Combinatorial} \label{ss:CD}
All of our examples where a path satisfies Equation \eqref{eq:hyp-cond1} but does not correspond to a stable irreducible component have the properties that
\begin{enumerate}

\item For some $k$, the point $(a_1+\cdots+a_{2k-1,}a_2+\cdots+a_{2k})$ is close to the diagonal (i.e. the line from $(0,0)$ to $(n,m)$), and

\item Immediately after that point there is some very unusual behavior. 
\end{enumerate}
It is well known that a large random rational Dyck path is usually far from the diagonal. For example, the following is immediate from \cite[Theorem 7.1]{BrMa}:
\begin{Proposition} \label{prop:2r}
For a random rational Dyck path from $(0,0)$ to $(k+1,k)$, the expected number of times the path visits a point $(a,b)$ with $b-a=r$ approaches $4r+4$ as $k$ approaches infinity. \qed
\end{Proposition}

\noindent So the expected number of visits to a given distance $r$ from the diagonal is bounded independent of $k$, and is linear in $r$. To get error bounds on our estimates, one would need to show that the probability of observing unusual enough behavior to cause a violation of stability near a point $(a,b)$ with $b-a=r$ decreases fast enough with $r$ so that the sum of the errors stays small. Another issue is that Proposition \ref{prop:2r} is for all Dyck paths, and restricting to those satisfying \eqref{eq:hyp-cond1} will have some effect.

\subsection{Representation theoretic} \label{ss:RD}
If $x, y$ is such that $\frac{y}{x} > \frac{\sqrt{r^2-4}+r}{2}$, there is a maximal finite irreducible dimensional $\Lambda$-module $T$ such that
\begin{itemize}
\item The 1-head of $T$ has dimension $y$. 

\item Let $U$ be the submodule of $T$ obtained by removing the $1$ head. Then the $0$-head of $U$ has dimension $x$. 
\end{itemize}
For example if $x=14$, $y=37$, this is 
$$ 0^11^30^81^{21} \oplus 1 0^31^8 \oplus 1 0^31^8,$$
where each summand is described via its socle filtration. 
There is of course a similar statement with the roles of the $1$-head and the $0$-head interchanged. 
One interpretation of Theorem \ref{th:good-cond} is that a pair $a_k=x, a_{k+1}=y$ with $\frac{y}{x} > \frac{\sqrt{r^2-4}+r}{2}$ forces the existence of a sub-module which is a quotient of this module, which forces a violation of stability. 

In other cases there can be infinite modules satisfying these conditions. For instance, consider $x=5, y=13$. There is a module given by 
\begin{equation} \label{eq:bgm}
\cdots0^{34} 1^{13} 0^51^20^11^10^21^50^{13},
\end{equation}
where the notation now means the $0$-head is $0^{13}$, the submodule not including that has head $1^5$, and so on. If $T$ has string data with $a_k=5, a_{k+1}=13$, and $k$ odd, some quotient of this must still be a submodule of $T$. But this quotient need not violate stability. Theorem \ref{th:cond2} can be interpreted as giving conditions when a quotient must in fact violate stability. But this can only happen when some segment of the path both starts very close to the diagonal and has weight close to the boundary of the imaginary cone, and seeing both of these together should be rare. 

The module in \eqref{eq:bgm} explains Example \ref{ex:type3}: Any quotient of the above that kills $\cdots 0^{34} 1^{13} 0^51^20^1$ must violate stability. Thus no data with $a_3=5, a_4=13$ can correspond to a stable component. For this sort of things to happen $\frac{y}{x} $ must be close to $ \frac{\sqrt{r^2-4}+r}{2}$.

\section*{Appendix 1: Computational evidence. With Colin Williams.} \label{sec:data}

We wrote Python code \cite{code} to calculate our upper bounds by counting Dyck paths satisfying Theorems \ref{th:good-cond} and \ref{th:cond2}, and compared to known root multiplicities from \cite[Figure 2]{FN} and \cite[Chapter 11]{Kac:1990}. 
For roots of the form $(n+1) \alpha_0+n \alpha_1$, the bound using Theorem \ref{th:good-cond} is exact up to $n=6$, and the bound using Theorem \ref{th:cond2} is exact up to $n=14$. For $16\alpha_0+15 \alpha_1$ the actual multiplicity is 815214, and the two bounds are 837218 
(over by 22004 or 2.7\%) and 815215 (over by 1). 
Our bounds are not as tight for  $15\alpha_0+16 \alpha_1$: the multiplicity is still 815214, but our bounds are 1234431 (over by 
419217 or 51.4\%)  
and 817505 (over by 2291 or $0.3\%$). 
The method is also successful on other roots. For instance, for $\beta=15\alpha_0+11\alpha_1,$ the two estimates are 
23868 and 23750, and 23750 is correct. 

We also used Monte-Carlo methods to consider larger roots, in particular $\beta=51 \alpha_0+50 \alpha_1$ and  $\beta=50 \alpha_0+51 \alpha_1$ . The actual multiplicity, calculated using Peterson's formula, implemented in sage by Judge \cite{sageticket18000}, is $2.03935 \times 10^{23}$. This is about a third of Avogadro's number, way too big to count. Instead we randomly sampled to estimate the fraction of Dyck paths satisfying the conditions. Multiplying by the number of rational Dyck paths gives an estimate of each upper bound. Here is the result of the largest samples we used (each took about 24 hours on a 2018 laptop):
$$\begin{array}{|c|c|c|c|c|c|}
\hline
\text{Root}
& \begin{array}{cc} \text{Paths} \\ \text{sampled} \end{array} 
& \begin{array}{cc} \text{Satisfied} \\ \text{Th \ref{th:good-cond} } \end{array} 
 & \begin{array}{cc}\text{Th \ref{th:good-cond}} \\  \text{Estimate}  \end{array} 
& \begin{array}{cc} \text{Satisfied} \\ \text{Th \ref{th:cond2}} \end{array} 
 & \begin{array}{cc} \text{Th \ref{th:cond2}}  \\  \text{Estimate} \end{array}\\
\hline
51\alpha_0+50\alpha_1&10^9 & 112637 &2.2283 \times 10^{23} & 103219 & 2.0419 \times 10^{23} \\
\hline
50\alpha_0+51\alpha_1 & 10^9 & 171935 &3.4013 \times 10^{23} & 103504 & 2.0476 \times 10^{23} \\
\hline
\end{array}
$$
For $51\alpha_0+50\alpha_1$, the first estimate is over by $9.2\%$ and the second by $0.13\%.$
The number of paths satisfying either theorem is roughly a Poisson random variable with standard deviation about $0.32\%$. Thus was can say with high confidence that our bound from Theorem \ref{th:good-cond} is overestimating the multiplicity by between 8.5\% and 10\%, and that the bound from Theorem \ref{th:cond2} is correct to within 1\%. 

For $50\alpha_0+51\alpha_1$, the estimates are off by $67\%$ and $0.4\%$ respectively ($\pm$ about 0.64\%). Perhaps this is even stronger evidence that the bounds stay quite good, since for roots of the form $n \alpha_0+(n+1) \alpha_1$ there is already non-trivial error at $n=15$, and this has not grown much by $n=50$.

\bibliographystyle{amsplain-ac}
\bibliography{mybib.bib}

\end{document}